\documentclass[reqno,12pt]{amsart}
\usepackage{amsmath,amsthm,amssymb,amsfonts,amscd}
\newtheorem{thm}{Theorem}[section]
\newtheorem{defn}[thm]{Definition}
\newtheorem{exmp}[thm]{Example}
\newtheorem{exdefn}[thm]{Example-Definition}
\newtheorem{cor}[thm]{Corollary}
\newtheorem{lem}[thm]{Lemma}
\newtheorem{prop}[thm]{Proposition}
\newtheorem{conj}[thm]{Conjecture}

\newtheorem{prob}[thm]{Problem}
\newtheorem{rem}[thm]{Remark}

\newtheorem*{pf}{Proof}
\numberwithin{equation}{section}
\def\C{{\mathbb C}}

\def\Q{{\mathbb Q}}
\def\R{{\mathbb R}}
\def\Z{{\mathbb Z}}
\def\A{{\mathcal A}}

\def\p{{\partial }}

\begin{document}
\title{Matrix Factorizations and Representations of Quivers I}
\author{Atsushi Takahashi}
\dedicatory{Dedicated to Professor Kyoji Saito on the occasion of his 
60th birthday}
\address{Research Institute for Mathematical Sciences, Kyoto University,
Kyoto 606-8502, Japan}
\email{atsushi@kurims.kyoto-u.ac.jp}
\begin{abstract}
This paper introduces a mathematical definition of the category of 
D-branes in Landau-Ginzburg orbifolds in terms of $A_\infty$-categories.
Our categories coincide with the categories of (graded)
matrix factorizations for quasi-homogeneous polynomials.
After setting up the necessary definitions, 
we prove that our category for the polynomial $x^{n+1}$ is 
equivalent to the derived category of 
representations of the Dynkin quiver of type $A_{n}$. 
We also construct a special stability condition for the triangulated category
in the sense of T.~Bridgeland, which should be the "origin" of the 
space of stability conditions.
\end{abstract}
\maketitle
%
%
\section{Introduction}
This paper introduces new triangulated categories associated to 
quasi-homogeneous polynomials which define isolated singularities 
only at the origin and relates those categories 
with the derived categories of representations of quivers.
Our motivation comes from K.~Saito's theory of primitive forms, 
especially from a problem in his study on regular weight systems 
and generalized root systems \cite{sa:1}.
We will explain the problem below.
Let $(a,b,c;h)$ be a quadruple of positive integers such that the 
function 
$$
\chi(T):=\frac{(T^{h-a}-1)(T^{h-b}-1)(T^{h-c}-1)}{(T^a-1)(T^b-1)(T^c-1)}
$$
has no poles.
Such a quadruple $W:=(a,b,c;h)$ is called a regular weight system.
It is known that $W$ is a regular weight system if and only if 
we have at least one polynomial $f(x,y,z)\in\C[x,y,z]$ such that 
$$
ax\frac{\p f}{\p x}+by\frac{\p f}{\p y}+cz\frac{\p f}{\p z}=hf
$$
and 
$$
X_0:=\{(x,y,z)\in\C^3~|~f(x,y,z)=0\}
$$
has an isolated singularity only at the origin.
Note that the (restricted) map 
$$
f:\C^3\backslash f^{-1}(0) \to \C\backslash\{0\}
$$
is a topologically locally trivial fiber bundle, 
and the general fiber $X_1:=f^{-1}(1)$ (called the 
Milnor fiber) is an open 2-dimensional 
complex manifold whose second homology group $H_2(X_1,\Z)$ 
is a free $\Z$-module of rank $\mu:=(h-a)(h-b)(h-c)/abc=
\lim_{t\to 1}\chi(T)$.
Since $X_1$ is real 4-dimensional, $H_2(X_1,\Z)$ has an intersection
form
$$
I_{H_2(X_1,\Z)}:H_2(X_1,\Z)\times H_2(X_1,\Z)\to \Z.
$$
If $f$ is a defining polynomial of a simple (ADE) singularity, 
then $(H_2(X_1,\Z),-I_{H_2(X_1,\Z)})$ gives the root lattice 
of the finite root system corresponding to the singularity.
In \cite{sa:1}, it is shown that $(H_2(X_1,\Z),-I_{H_2(X_1,\Z)})$ with 
the set of vanishing cycles (which corresponds to the set of roots)
and the Milnor monodromy (which corresponds to the Coxeter transformation)
satisfies the axioms of the generalized root system which naturally 
extends the classical (finite) root systems.
Since both weight systems and generalized root systems are 
combinatorial, it is natural to propose the following 
problem.
\begin{prob}$($\cite{sa:1}$)$\\
Construct directly from a regular weight system $W$,
without passing through the homology group $H_2(X_1,\Z)$ of the 
Milnor fiber, arithmetically or combinatorially, the generalized 
root system of the vanishing cycles.
\end{prob}
The purpose of this paper is to develop a necessary tools 
in terms of $A_\infty$-categories and to give 
a partial answer to the above problem.
Let $k$ be a field of characteristic zero.
First, we introduce a notion of $\Q$-graded $A_\infty$-categories 
over $k$ (Definition \ref{defn:qainfty}) 
in order to consider the polynomial ring $\C[x_1,\dots x_n]$ with 
a polynomial $f$ satisfying the quasi-homogeneous condition
\begin{equation}
\sum_{i=1}^n \frac{2a_i}{h}\cdot x_i\frac{\p f}{\p x_i} =2f,
\end{equation}
where $a_1,\dots, a_n$ and $h$ are positive integers such that 
the greatest common divisor of them is $1$, 
as a usual $\Z_2$-graded $A_\infty$-category with $m_0(1)=f$ 
and an "extra $\frac{2}{h}\cdot\Z$-grading".
We shall denote the $\Q$-graded $A_\infty$-category 
defined by $f\in\C[x_1,\dots,x_n]$ by $\A_f$ 
(Example-Definition \ref{defn:af}).
Next, we consider the category of twisted complexes over $\Q$-graded 
$A_\infty$-categories (Proposition \ref{prop:tw}) 
and the derived category of $\Q$-graded $A_\infty$-categories
(Definition \ref{defn:derived}).
The important fact is that the twisted complexes over $\A_f$ coincide 
with matrix factorizations of $f$ introduced by Eisenbud \cite{e:1} 
in his study of maximal Cohen-Macaulay modules.
Since we consider the quasi-homogeneous polynomial $f$, 
we have a group action ($\Z$-action) on the category 
of matrix factorizations.
Inspired by the work by Hori and Walcher \cite{hw:1}, 
we introduce the $\Z$-equivariant derived category of $\A_f$ 
denoted by $D^b_{\Z}(\A_f)$ (Definition \ref{defn:zderived}).
We can now propose a conjecture to K.~Saito's problem.
\begin{conj}
Let $W$ be a regular weight system and $f$ be a 
quasi-homogeneous polynomial attached to $W$.
Assume $W$ has a dual regular weight system $W^*=(a^*,b^*,c^*;h)$ 
in the sense of \cite{sa:2} and 
let $f^*$ be a quasi-homogeneous polynomial attached to $W^*$. 
Then the following should hold.
\begin{enumerate}
\item $D^b_{\Z}(\A_{f^*})$ is generated as a triangulated category 
by objects $\{E_1,\dots,E_{\mu}\}$ such that 
\begin{equation}
{\rm Hom}_{D^b_{\Z}(\A_{f^*})}(E_i,E_j)=0,\quad \text{if }i>j,\quad 
{\rm Hom}_{D^b_{\Z}(\A_{f^*})}(E_i,E_j[k])=0,\quad k\ne 0,{}^\forall i,j.
\end{equation}
That is to say, $D^b_{\Z}(\A_{f^*})$ is generated by a strongly exceptional 
collection.
\item $D^b_{\Z}(\A_{f^*})$ has the Serre functor $S$ such that 
$S^h\simeq [3h-2a^*-2b^*-2c^*]$ where $[1]$ is the shift functor on 
$D^b_{\Z}(\A_{f^*})$.
\item Let $a_{ij}:=\chi(E_i,E_j)=\dim_k
{\rm Hom}_{D^b_{\Z}(\A_{f^*})}(E_i,E_j)$.
Put $A:=(a_{ij})$ and $I_{K_0(D^b_{\Z}(\A_{f^*}))}=A^{-1}+{}^tA^{-1}$.
Then $(K_0(D^b_{\Z}(\A_{f^*})),I_{K_0(D^b_{\Z}(\A_{f^*}))})$ 
is isomorphic to $(H_2(X_1,\Z),-I_{H_2(X_1,\Z)})$ as a lattice.
\end{enumerate}
\end{conj}
This conjecture is based on the relation between the duality of 
regular weight systems and the mirror symmetry of 
Landau-Ginzburg orbifolds (see \cite{t:1}).
We do not discuss this background in detail here but 
we write the following diagram for reader's convenience.
\begin{equation*}
\begin{CD}
\text{Quasi-homogeneous polynomial }f \text{ for }W &
\stackrel{Milnor fiber}{\Longrightarrow} &
\{\text{Vanishing cycles in }X_1=f^{-1}(1)\}\\
\text{Duality of weights + Orbifold}\downarrow & & ||\\
\{\text{B-branes in LG orbifold }W^*//(\Z/h\Z)\} &
\stackrel{Mirror Symmetry}{\simeq} &
\{\text{A-branes in LG model for }W\}
\end{CD}
\end{equation*}
For ADE singularities, we know that $W\simeq W^*$ 
and the generalized root systems for them are 
the classical finite root systems.
Therefore, we may expect that the following conjecture should hold.
\begin{conj}
Let $W$ be a regular weight system corresponding to an 
ADE singularity
and $f$ be a quasi-homogeneous polynomial attached to $W$.
Then $D^b_{\Z}(\A_f)$ is equivalent as a triangulated category 
to the bounded derived category of finite dimensional 
representations of Dynkin quiver corresponding to the type of 
singularity of $f$.
\end{conj}
\noindent
This is also expected from the homological mirror symmetry phenomena
for ADE singularities studied by Seidel. 
See \cite{se:2} for details.
\bigskip
In this paper, we will prove the conjecture for $A_n$-singularities 
(Theorem \ref{thm:quiver}), 
where we reduce to the case $f:=x^{n+1}\in \C[x]$ 
by Kn\"{o}rrer's periodicity \cite{k:1} (see also \cite{o:1}). 
We will give a proof of the above conjecture for general cases 
in a separate paper \cite{kst:1}.
Finally, we will construct a special stability condition 
for the triangulated category $D^b_{\Z}(\A_f)$ in the sense of 
T.~Bridgeland \cite{b:1} for $f=x^h$.
We can naturally introduce in our formulation 
the phase of objects (Definition \ref{defn:phase}) 
and the central charge $Z_\omega$ (Definition \ref{defn:centralcharge}).
\bigskip
While our preparation of this paper, 
two papers related to our work appeared.
One is the paper \cite{w:1} by J.~Walcher where he studies 
from physical point of view the similar categories and 
the stability conditions on them (his notion of "R-stability").
Another is the paper \cite{o:2} by D.~Orlov where he studies 
the triangulated category for singularities with a $\C^*$-action,
which is equivalent to our category $D^b_\Z(\A_f)$.
\bigskip
\noindent
{\bf Acknowledgement}\\
\indent
I am grateful to Kentaro Hori, Hiroshige Kajiura, Masaki Kashiwara and Kyoji 
Saito for very useful discussions.
This work was partly supported by Grant-in Aid for Scientific Research 
grant numbers 14740042 of the Ministry of Education, Science and Culture 
in Japan.
\section{$\Q$-graded $A_\infty$-categories}
In this section, we set up several definitions which we will 
use in the later sections.
Let $k$ be a field of characteristic zero.
\begin{defn}\label{defn:qainfty}
Let $h$ be a positive number. 
A $\Q$-graded $A_\infty$-category $\A$ of index $h$ is a collection of 
the following data.
\begin{enumerate}
\item A set of objects $Ob(\A)$, 
\item A set of homomorphisms, a $\Q\times \Z_2$ graded 
$k$-linear vector space for each $a,b\in Ob(\A)$
$$
\A(a,b)=\bigoplus_{q\in\Q}\A^q(a,b)_+\oplus\A^q(a,b)_-,
$$
such that 
$$
\A^q(a,b)_+=0,\quad q\notin\frac{2}{h}\Z,\quad 
\A^q(a,b)_-=0,\quad q-1\notin\frac{2}{h}\Z.
$$
We call the subspaces
$$
\A(a,b)_+:=\bigoplus_{q\in\frac{2}{h}\Z}\A^q(a,b)_+,
\quad\A(a,b)_-:=\bigoplus_{q-1\in\frac{2}{h}\Z}\A^q(a,b)_-
$$ 
the even and the odd subspaces.
\item for $n\ge 0$, $k$-multilinear maps
$$
m_n^\A:\A(a_{n-1},a_n)\otimes\dots\otimes\A(a_0,a_1)\to \A(a_0,a_n),\quad 
a_i\in Ob(\A),
$$
of degree $2-n$ with respect to the $\Q$-grading which is even $($odd$)$ 
with respect to the $\Z_2$-grading when $n$ is even $($odd$)$, 
where $m_0^\A$ is a map 
$$
m_0^\A:k\to \A(a,a).
$$
The multilinear maps satisfy the following $(A_\infty$-relation$)$. 
For fixed $n$, we have 
\begin{multline}\label{eq:A-infinity}
\sum_{r+s+t=n}\sum_{r+1+t=u}
(-1)^{|x_{1}|+\dots+|x_{r}|+r}m_u^\A(x_{r+s+t}\otimes\dots\otimes x_{r+s+1}
\otimes \\
m_s^\A(x_{r+s}\otimes\dots\otimes x_{r+1})\otimes
x_{r}\otimes\dots\otimes x_{1})=0,
\end{multline}
where $|x_i|$ is the parity of the morphism defined by
\begin{equation}
|x_i|:=
\left\{\begin{array}{l}
0,\quad x_i\in \A(a_{i-1},a_i)_+,\\
1,\quad x_i\in \A(a_{i-1},a_i)_- 
\end{array}.\right.
\end{equation}
\end{enumerate}
\end{defn}
\begin{rem}
$\Q$-graded $A_\infty$-category of index $1$ is nothing 
but an $A_\infty$-category with the usual $\Z$-grading, 
we call it a $\Z$-graded $\A_\infty$-category or 
simply an $\A_\infty$-category.
See \cite{f:1} and \cite{se:1} for details of homological algebra of 
$\A_\infty$-categories.
\end{rem}
We write down explicitly the relation \eqref{eq:A-infinity} 
when $m_n^\A=0$ for $n\ge 3$.
For $x,y,z\in\oplus_{a,b}\A(a,b)$, we have
$$
m_1^\A(m_0^\A(1))=0,
$$
$$
m_1^\A(m_1^\A(x))=(-1)^{|x|}m_2^\A(m_0^\A(1)\otimes x)
-m_2^\A(x\otimes m_0^\A(1)),
$$
$$
m_1^\A(m_2^\A(x\otimes y))=
(-1)^{|y|}m_2^\A(m_1^\A(x)\otimes y)-m_2^\A(x\otimes m_1^\A(y)),
$$
$$
m_2^\A(m_2^\A(x\otimes y)\otimes z)=
(-1)^{|z|}m_2^\A(x\otimes m_2^\A(y\otimes z)).
$$
Put $u:=-m_0^\A(1)$, $d(x):=(-1)^{|x|+1}m_1^\A(x)$ and 
$x\cdot y:=(-1)^{|y|}m_2^\A(x\otimes y)$. 
Then a triple $(u,d,\cdot )$ defines on $\oplus_{a,b}\A(a,b)$ 
a curved differential graded (CDG) algebra structure \cite{kl:1}.
\begin{exdefn}\label{defn:af}
Let $f\in \C[x_1,\dots x_n]$ be a polynomial which satisfies the 
following quasi-homogeneous condition$:$
\begin{equation}
\sum_{i=1}^n \frac{2a_i}{h}\cdot x_i\frac{\p f}{\p x_i} =2f,
\end{equation}
where $a_1,\dots, a_n$ and $h$ are positive integers such that 
the greatest common divisor of them is $1$.
We denote by $\A_f$ the $\Q$-graded $A_\infty$-category of index $h$ 
defined as follows: 
$$
Ob(\A_f)=\{a\},
$$
$$
\A_f(a,a):=\C[x_1,\dots x_n],\quad \A_f(a,a)_{-}:=0,
$$
$$
m_0^{\A_f}(1):=f\in \A^2_f(e,e)_+,
\quad m_1^{\A_f}:=0,
$$
$$
m_2^{\A_f}(\alpha\otimes\beta):=\alpha\cdot\beta,\quad \alpha,\beta\in
\C[x_1,\dots x_n],
$$ 
where $\cdot$ is the usual product on $\C[x_1,\dots x_n]$.
\end{exdefn}
In the above example, we have the special element $1\in\A^0_f(a,a)_+$ 
which defines a unit of the algebra $\C[x_1,\dots x_n]$.
It is well-known that the notion of units in $A_\infty$-categories
can be introduced as follows:
\begin{defn}
Let $a\in Ob(\A)$. $e_a\in\A^0(a,a)_+$ is called a unit if 
\begin{equation}
m_2^\A(x,e_a)=x,\quad m_2^\A(e_a,y)=(-1)^{|y|}y,
\end{equation}
for $x\in\A(a,b)$ and $y\in\A(b,a)$, and for $n\ne 2$
\begin{equation}
m_n^\A(x_1,\dots,x_n)=0, 
\end{equation}
if one of $x_i$ coincides with $e_a$.
\end{defn}
\begin{rem}
It is easy to check that if a unit exists then it is unique.
\end{rem}
\begin{defn}
$\Q$-graded $A_\infty$-category is called unital if each object has a unit. 
\end{defn}
Let $\A$ be a $\Q$-graded $A_\infty$-category of index $h$.
We can construct another $\Q$-graded $A_\infty$-category $\overline{\A}$ 
of index $h$ from $\A$ as follows:
\begin{defn}
Let $\A$ be a unital $\Q$-graded $A_\infty$-category of index $h$.
\begin{enumerate}
\item a set of objects $Ob(\overline{\A})$ is given by
\begin{equation}
Ob(\overline{\A}):=\left\{a\{\frac{2k}{h}\},a\in Ob(\A),k\in\Z,\quad 
b\{\frac{2l}{h}\}[-1],b\in Ob(\A),l\in\Z\right\},
\end{equation}
\item a set of homomorphisms is given by 
\begin{equation}
\overline{\A}^q\left(a\{\frac{2k}{h}\},b\{\frac{2l}{h}\}\right)_\pm:=
\A^{q+\frac{2(l-k)}{h}}(a,b)_\pm,
\end{equation}
\begin{equation}
\overline{\A}^q\left(a\{\frac{2k}{h}\},b\{\frac{2l}{h}\}[-1]\right)_\pm:=
\A^{q+\frac{2(l-k)}{h}-1}(a,b)_\mp,
\end{equation}
\begin{equation}
\overline{\A}^q\left(a\{\frac{2k}{h}\}[-1],b\{\frac{2l}{h}\}\right)_\pm:=
\A^{q+\frac{2(l-k)}{h}+1}(a,b)_\mp,
\end{equation}
\begin{equation}
\overline{\A}^q\left(a\{\frac{2k}{h}\}[-1],b\{\frac{2l}{h}\}[-1]\right)_\pm:=
\A^{q+\frac{2(l-k)}{h}}(a,b)_\pm.
\end{equation}
\item $k$-multilinear maps $m_n^{\overline{\A}}$ 
are defined using those on $\A$ 
with additional signs as follows.
For $x_1\in\overline{\A}(a_0,a_1) ,\dots, x_n\in \overline{\A}(a_{n-1},a_n)$,
$$
m^{\overline{\A}}_n(x_n\otimes \dots\otimes x_1):=(-1)^{|a_0|}
m_n^{\A}(x_n\otimes \dots\otimes x_1),
$$
where we regard $x_i$ in the right hand side 
as a homomorphism in $\A$ by the above definition $(ii)$
and $|a_0|$ is the parity of $a_0$ defined by
\begin{equation}
|a_0|:=
\left\{\begin{array}{l}
0,\quad a_0=a\{\frac{2k}{h}\},\quad a\in Ob(\A), k\in\Z\\
1,\quad a_0=a\{\frac{2l}{h}\}[-1],\quad a\in Ob(\A), l\in\Z
\end{array}.\right.
\end{equation}
\end{enumerate}
\end{defn}
$A_\infty$-functors for $\Q$-graded $A_\infty$-categories 
can be defined in an obvious way.
There are the "translation" functor $\{\frac{2}{h}\}$ and 
the shift functor $[1]$ on $\overline{\A}$.
\begin{prop}
The following functors $\{\frac{2}{h}\}$ and $[1]$ define 
autoequivalences of $\overline{\A}:$
\begin{equation}
\{\frac{2}{h}\}\left(a\{\frac{2k}{h}\}\right):=a\{\frac{2(k+1)}{h}\},\quad 
\{\frac{2}{h}\}\left(a\{\frac{2k}{h}\}[-1]\right):=a\{\frac{2(k+1)}{h}\}[-1],
\end{equation}
and 
\begin{equation}
[1]\left(a\{\frac{2k}{h}\}\right):=a\{\frac{2(k+h)}{h}\}[-1],\quad 
[1]\left(a\{\frac{2k}{h}\}[-1]\right):=a\{\frac{2k}{h}\}.
\end{equation}
Put $\{\frac{2k}{h}\}:=\{\frac{2}{h}\}^k$ and $[l]:=[1]^l$ for $k,l\in\Z$.
We have the relation $\{\frac{2h}{h}\}=[2]$.
\qed
\end{prop}
Let $\overline{\A}$ be an $\Q$-graded $A_\infty$-category of index $h$.
Consider the $\Q$-graded $A_\infty$-category $\widetilde{\A}$ of index $h$ 
whose set of objects is the set of finite (formal) direct sums 
of objects of $\overline{\A}$, 
\begin{equation}
Ob(\widetilde{\A}):=\left\{
a=\bigoplus_i a_i\{\frac{2k_i}{h}\} \oplus
\bigoplus_j a_j\{\frac{2l_j}{h}\}[-1],
\quad a_i,a_j\in Ob(\overline{\A}),\quad k_i,l_j\in\Z\right\},
\end{equation}
whose set of homomorphisms is
\begin{multline}
\widetilde{\A}(a,b):=\bigoplus_{i_1,i_2}
\overline{\A}(a_{i_1}\{\frac{2k_{i_1}}{h}\},b_{i_2}\{\frac{2k_{i_2}}{h}\})
\oplus\bigoplus_{i_1,j_2}
\overline{\A}(a_{i_1}\{\frac{2k_{i_1}}{h}\},
b_{j_2}\{\frac{2l_{j_2}}{h}\}[-1])\\
\oplus\bigoplus_{j_1,i_2}
\overline{\A}(a_{j_1}\{\frac{2l_{j_1}}{h}\}[-1],b_{i_2}\{\frac{2k_{i_2}}{h}\})
\oplus\bigoplus_{j_1,j_2}
\overline{\A}(a_{j_1}\{\frac{2l_{j_1}}{h}\},
b_{j_2}\{\frac{2l_{j_2}}{h}\}[-1]),
\end{multline}
and whose $k$-linear maps are defined by those on $\overline{\A}$ 
using the natural "matrix multiplication" rule.
\begin{defn}\label{defn:twistcomp}
Take an object $a \in Ob(\widetilde{\A})$ and 
$Q\in\widetilde{\A}(a,a)_-$.
$(a;Q)$ is called a twisted complex if $Q$ satisfies 
the Maurer-Cartan equation
\begin{equation}\label{eq:mc}
\sum_{n\ge 0}m^{\widetilde{\A}}_n(Q^{\otimes n})=0.
\end{equation}
The set of all twisted complexes is denoted by 
$Ob(Tw(\A))$.
If $Q\in\widetilde{\A}^1(a,a)_-$ in addition, 
then $(a;Q)$ is called a graded twisted complex and we denote 
the set of all graded twisted complexes by $Ob(Tw_\Z(\A))$.
\end{defn}
An assumption is necessary for the equation \eqref{eq:mc} to make sense.
If $m_0^\A=0$, then it is usually introduced that 
the notion of one-sided twisted complexes which makes the sum in 
the equation \eqref{eq:mc} finite.
However we study in this paper the case when $m_0^\A\ne0$, 
we shall assume that our $A_\infty$-categories have no higher product, 
in other words, $m_n^\A=0$ for all $n\ge 3$.
We often write $Q$ in the following form:
\begin{equation}
Q=
\begin{pmatrix}
Q_{++} & Q_{-+}\\
Q_{+-} & Q_{--},
\end{pmatrix}
\end{equation}
where
\begin{equation}
Q_{\pm\pm}\in \widetilde{\A}(a_\pm,a_\pm)_-, 
Q_{\pm\mp}\in \widetilde{\A}(a_\pm,a_\mp)_+, 
\end{equation}
and $a_\pm$ are given by the following decomposition
\begin{equation}
a=a_+ + a_-[-1], \quad a_+=\bigoplus_{i}a_{+,i}\{\frac{2k_i}{h}\},
\quad a_-=\bigoplus_{i}a_{-,i}\{\frac{2l_i}{h}\}.
\end{equation}
\begin{rem}
If $\A$ is a unital $\Q$-graded $A_\infty$-category of index $h$ 
with $m_n^\A=0$, $n\ge 3$, then there exists at least one twisted complex 
for each object $a\in Ob(\widetilde{\A})$.
Indeed, 
\begin{equation}
Q|_{ij}=\left\{\begin{array}{l}
0, \quad \text{for }i,j\ \text{such that }a_{+,i}\ne a_{-,j} \text{ and }
a_{+,i}\{\frac{2\cdot h}{h}\}\ne a_{-,j}, \\
\begin{pmatrix}
0 & m_0^\A(1)\\
e_{a_{+,i}} & 0
\end{pmatrix},\quad \text{for }i,j\ \text{such that }a_{+,i}= a_{-,j},\\
\begin{pmatrix}
0 & e_{a_{+,i}}\\
m_0^\A(1) & 0
\end{pmatrix},\quad \text{for }i,j\ \text{such that }
a_{+,i}\{\frac{2\cdot h}{h}\}= a_{-,j},
\end{array} \right.
\end{equation}
is a twisted complex.
\end{rem}
\begin{exmp}
Since $\A_f$ has no odd homomorphisms, 
each twisted complex $(a=a_+\oplus a_-[-1];Q_a)$ has the following form 
$$
Q_a:=
\begin{pmatrix}
0 & Q_{-+}\\
Q_{+-} & 0
\end{pmatrix},
\quad Q_{+-}\in\widetilde{\A}(a_+,a_-)_+,\quad 
Q_{-+}\in\widetilde{\A}(a_-,a_+)_+.
$$
The Maurer-Cartan equation \eqref{eq:mc} becomes 
\begin{equation}
f\cdot{\rm Id} -Q_a^2=0.
\end{equation}
This is exactly the same equation which first studied by 
Eisenbud \cite{e:1} in his work on maximal Cohen-Macaulay modules. 
$Q_a$ is called a matrix factorization of $f$.
\end{exmp}
Let $\A$ be a unital $\Q$-graded $A_\infty$-category of index $h$ 
with $m_n^\A=0$, $n\ge 3$.
\begin{defn}\label{defn:twistmor}
Let $\alpha:=(a;Q_a)$ and $\beta:=(b;Q_b)$ be twisted complexes.
We first put 
\begin{equation}
Tw(\A)(\alpha,\beta):=\widetilde{\A}(a,b)_+\oplus 
\widetilde{\A}(a,b)_-.
\end{equation}
We define a $k$-multilinear maps $m_n^{Tw(\A)}(n=0,1,2)$ by 
\begin{equation}
m_0^{Tw(\overline{\A})}(1):=0,
\end{equation}
\begin{equation}
m_1^{Tw(\A)}(\Phi) :=m_1^{\widetilde{\A}}(\Phi)
+m_2^{\widetilde{\A}}(Q_b\otimes\Phi)
+m_2^{\widetilde{\A}}(\Phi\otimes Q_a),
\end{equation}
where $\Phi\in Tw(\A)(\alpha,\beta)$ and 
\begin{equation}
m_2^{Tw(\A)}(\Psi_2\otimes \Psi_1):=m_2^{\widetilde{\A}}
(\Psi_2\otimes \Psi_1),
\end{equation}
for $\Psi_1\in Tw(\A)(\alpha_0,\alpha_1)=\widetilde{\A}(a_0,a_1)$ 
and $\Psi_2\in Tw(\A)(\alpha_1,\alpha_2)=\widetilde{\A}(a_1,a_2)$.
\end{defn}
We often write the spaces of morphisms in the matrix form:
$$
Tw(\A)(\alpha,\beta)_\pm=
\begin{pmatrix}
\widetilde{\A}(a_+,b_+)_\pm & 
\widetilde{\A}(a_-,b_+)_\mp\\
\widetilde{\A}(a_+,b_-)_\mp & 
\widetilde{\A}(a_-,b_-)_\pm
\end{pmatrix}.
$$
\begin{lem}
$(m_1^{Tw(\A)})^2=0$.
\end{lem}
\begin{pf}
For $\Phi_\pm\in Tw(\A)(\alpha,\beta)_\pm$, we have
\begin{align*}
(m_1^{Tw(\A)})^2(\Phi_\pm) =&
(m_1^{\widetilde{\A}})^2(\Phi_\pm)
+m_1^{\widetilde{\A}}(m_2^{\widetilde{\A}}(Q_b\otimes\Phi_\pm))
+m_1^{\widetilde{\A}}(m_2^{\widetilde{\A}}(\Phi_\pm\otimes Q_a))\\
&+ m_2^{\widetilde{\A}}(Q_b\otimes (m_1^{\widetilde{\A}}(\Phi_\pm) 
+m_2^{\widetilde{\A}}(Q_b\otimes\Phi_\pm)
+ m_2^{\widetilde{\A}}(\Phi_\pm\otimes Q_a)))\\
&+  m_2^{\widetilde{\A}}((m_1^{\widetilde{\A}}(\Phi_\pm) 
+m_2^{\widetilde{\A}}(Q_b\otimes\Phi_\pm)
+ m_2^{\widetilde{\A}}(\Phi_\pm\otimes Q_a))\otimes Q_a)\\
=&\pm m_2^{\widetilde{\A}}(m^{\widetilde{\A}}_0(1)\otimes\Phi_\pm )
-m_2^{\widetilde{\A}}(\Phi_\pm\otimes m^{\widetilde{\A}}_0(1) )
\pm m_2^{\widetilde{\A}}(m_1^{\widetilde{\A}}(Q_b)\otimes \Phi_\pm)\\
& - m_2^{\widetilde{\A}}(\Phi_\pm\otimes m_1^{\widetilde{\A}}(Q_a))
\pm m_2^{\widetilde{\A}}(m_2^{\widetilde{\A}}(Q_b^{\otimes 2})\otimes \Phi_\pm)
- m_2^{\widetilde{\A}}(\Phi_\pm\otimes m_2^{\widetilde{\A}}(Q_a^{\otimes 2}))\\
=& \pm m_2^{\widetilde{\A}}((m_0^{\widetilde{\A}}(1)+m_1^{\widetilde{\A}}
(Q_b)+m_2^{\widetilde{\A}}(Q_b^{\otimes 2})\otimes \Phi_\pm)\\
&- m_2^{\widetilde{\A}}(\Phi_\pm\otimes(m_0^{\widetilde{\A}}(1
+m_1^{\widetilde{\A}}(Q_a)+m_2^{\widetilde{\A}}(Q_a^{\otimes 2}) )\\
=&0.
\end{align*}
\qed
\end{pf}
\begin{lem}
For $\Phi\in Tw(\A)(\alpha_0,\alpha_1)$ and 
$\Psi\in Tw(\A)(\alpha_1,\alpha_2)$, we have
\begin{equation}
m_1^{Tw(\A)}(m_2^{Tw(\A)}(\Psi\otimes \Phi))=
(-1)^{|\Phi|}m_2^{Tw(\A)}(m_1^{Tw(\A)}(\Psi)\otimes \Phi)
-m_2^{Tw(\A)}(\Psi\otimes m_1^{Tw(\A)}(\Phi)).
\end{equation}
\end{lem}
\begin{pf}
\begin{align*}
&m_1^{Tw(\A)}(m_2^{Tw(\A)}(\Psi\otimes \Phi))\\
=&
m_1^{\widetilde{\A}}(m_2(\Psi\otimes \Phi))
+m_2^{\widetilde{\A}}(Q_{a_2}\otimes m_2^{\widetilde{\A}}(\Psi\otimes \Phi))
+m_2^{\widetilde{\A}}(m_2^{\widetilde{\A}}(\Psi\otimes \Phi)\otimes Q_{a_0})\\
=&(-1)^{|\Phi|}m_2^{\widetilde{\A}}(m_1^{\widetilde{\A}}(\Psi)\otimes \Phi))
-m_2^{\widetilde{\A}}(\Psi\otimes m_1^{\widetilde{\A}}(\Phi))\\
&(-1)^{|\Phi|}m_2^{\widetilde{\A}}
(m_2^{\widetilde{\A}}(Q_{a_2}\otimes\Psi)\otimes \Phi)
-m_2^{\widetilde{\A}}(\Psi\otimes m_2^{\widetilde{\A}}(\Phi\otimes Q_{a_0}))\\
=&(-1)^{|\Phi|}m_2^{Tw(\A)}(m_1^{Tw(\A)}(\Psi)\otimes \Phi)
-(-1)^{|\Phi|}m_2^{\widetilde{\A}}(m_2^{\widetilde{\A}}
(\Psi\otimes Q_{a_1})\otimes \Phi)\\
&-m_2^{Tw(\A)}(\Psi\otimes m_1^{Tw(\A)}(\Phi))
+m_2^{\widetilde{\A}}(\Psi\otimes 
m_2^{\widetilde{\A}}(Q_{a_1}\otimes\Phi))\\
=&(-1)^{|\Phi|}m_2^{Tw(\A)}(m_1^{Tw(\A)}(\Psi)\otimes \Phi)
-m_2^{Tw(\A)}(\Psi\otimes m_1^{Tw(\A)}(\Phi)).
\end{align*}
\qed\end{pf}
By the above two Lemmas, we have the following.
\begin{prop}\label{prop:tw}
Let $\A$ be a unital $\Q$-graded $A_\infty$-category of index $h$ 
with $m_n^\A=0$, $n\ge 3$.
A collection $Ob(Tw(\A))$, 
$Tw(\A)(\alpha,\beta)$ and 
$(m_0^{Tw(\A)},m_1^{Tw(\A)},m_2^{Tw(\A)})$
given by Definition \ref{defn:twistcomp} and Definition \ref{defn:twistmor} 
determines a structure of a differential $(\Z_2$-$)$graded category. 
We denote it by $Tw(\A)$.
\qed
\end{prop}
\begin{rem}
Note that the condition that $\A$ is $\Q$-graded is not necessary for 
the above definition of the category $Tw(\A)$ 
and the category $D^b(\A)$ below. We need only the $\Z_2$-grading. 
\end{rem}
\begin{defn}\label{defn:derived}
Let $\A$ be a unital $\Q$-graded $A_\infty$-category of index $h$ 
with $m_n^\A=0$, $n\ge 3$.
We construct the category $D^b(\A)$ called 
the bounded derived category of $\A$ as follows.
The set of objects is given by 
\begin{equation}
Ob(D^b(\A)):=Ob(Tw(\A)),
\end{equation}
and the set of homomorphisms is given by
\begin{align}
{\rm Hom}_{D^b(\A)}(\alpha,\beta):=&
{\rm Ker}(m_1^{Tw(\A)}:Tw(\A)(\alpha,\beta)_+\to 
Tw(\A)(\alpha,\beta)_-)\\
&\left/ {\rm Im}(m_1^{Tw(\A)}:Tw(\A)(\alpha,\beta)_-\to 
Tw(\A)(\alpha,\beta)_+)\right. .
\end{align}
\end{defn}
Let $T\in Tw(\A)(\alpha,\beta)_+$ be a $m_1^{Tw(\A)}$-closed 
homomorphism.
We define a mapping cone $C(T)$ as an object
\begin{equation}
C(T):=(a[1]\oplus b; Q_{C(T)}),\quad
Q_{C(T)}:=
\begin{pmatrix}
Q_{a[1]} & 0\\
T & Q_b
\end{pmatrix}
.
\end{equation}
$C(T)$ is well-defined since the Maurer-Cartan equation \eqref{eq:mc} 
for $Q_{C(T)}$ is equivalent to the equation 
$m_1^{Tw(\A)}(T)=0$ and 
the Maurer-Cartan equation \eqref{eq:mc} for $Q_a$ and $Q_b$. 
Note also that there are natural closed morphisms 
$$
\beta\to C(T),\quad C(T)\to \alpha.
$$
We define an exact triangle in the category $D^b(\A)$
as a triangle of the form 
\begin{equation}
\alpha\stackrel{T}{\to}\beta \to C(T)\to \alpha[1],
\end{equation}
for some $T\in {\rm Ker}(m_1^{Tw(\A)}:Tw(\A)(\alpha,\beta)_+\to 
Tw(\A)(\alpha,\beta)_-)$. 
\begin{thm}
The category $D^b(\A)$ endowed with a shift functor 
$[1]$ and the class of exact triangles defined above
becomes a triangulated category.
\end{thm}
\begin{pf}
The proof is essentiality the same as the known results 
in the usual situation.
See for example, \cite{bk:1},\cite{gm:1}, \cite{ks:1}, \cite{o:1}
and \cite{se:1}.
\qed\end{pf}
\begin{rem}
The twice of the shift functor $[2]:=[1]^{\otimes 2}$ is isomorphic to the identity functor 
in $D^b(\A)$.
\end{rem}
We shall add more objects to $D^b(\A)$ following \cite{se:1}.
\begin{defn}
Consider the category $D^\pi(\A)$ whose objects are pairs 
$(X,p)$ where $X\in Ob(D^b(\A))$ and $p\in {\rm Hom}_{D^b(\A)}(X,X)$
an idempotent endomorphism, and whose spaces of homomorphisms are
${\rm Hom}_{D^\pi(\A)}((X_0,p_0),(X_1,p_1)):=
p_1{\rm Hom}_{D^b(\A)}(X_0,X_1)p_0$.
The category $D^\pi(\A)$ is called the split-closed derived category of $\A$.
\end{defn}
It is known that $D^\pi(\A)$ is again a triangulated category 
(see \cite{bs:1}).
\begin{rem}
Since any projective module over the polynomial ring $\C[x_1,\dots,x_n]$ 
is a free module, we have $D^b(\A_f)\simeq D^\pi(\A_f)$.
Therefore, we shall study in this paper only the category $D^b(\A_f)$.
It is not difficult to see that our category $D^b(\A_f)$ is equivalent to 
the category of matrix factorizations, or equivalently, the category of 
maximal Cohen-Macaulay modules over $\C[x_1,\dots,x_n]/(f)$ 
without free summands studied by Eisenbud \cite{e:1}, Kn\"{o}rrer \cite{k:1},
Buchweitz-Greuel-Schreyer \cite{bgs:1}, Orlov\cite{o:1} 
and other people $($see the book by Yoshino \cite{y:1} for 
the details on maximal Cohen-Macaulay modules$)$.
Indeed, we can construct a functor from the category of matrix factorizations 
to $D^b(\A_f)$ once we choose a basis of the free module over 
the polynomial ring.
Note that any object isomorphic to a direct sum of the following objects
$$
\begin{pmatrix}
0 & f\\
1 & 0
\end{pmatrix},
\quad
\begin{pmatrix}
0 & 1\\
f & 0
\end{pmatrix}
$$
becomes the zero object of the category $D^b(\A_f)$.
\end{rem}
Next, we shall define the $\Z$-equivariant bounded derived category 
$D^b_{\Z}(\A)$ of $\A$.
Let $\alpha:=(a;Q_a)$ and $\beta:=(b;Q_b)$ be graded twisted complexes.
We put 
$$
Tw_{\Z}(\A)(\alpha,\beta):=\bigoplus_{q\in\Z}Tw_{\Z}^q(\A)(\alpha,\beta), 
$$
where 
\begin{equation}
Tw_{\Z}^q(\A)(\alpha,\beta):=\left\{\begin{array}{l}
\widetilde{\A}^q(a,b)_+,\quad q\in 2\Z,\\
\widetilde{\A}^q(a,b)_-,\quad q-1\in 2\Z. 
\end{array}\right.
\end{equation}
Since $Tw_{\Z}^q(\A)(\alpha,\beta)\subset Tw(\A)(\alpha,\beta)$, 
we can define a $k$-multilinear maps $m_n^{Tw_{\Z}(\A)}$ by restricting 
$m_n^{Tw(\A)}$ to the subspaces.
\begin{prop}
Let $\A$ be a unital $\Q$-graded $A_\infty$-category of index $h$ 
with $m_n^\A=0$, $n\ge 3$.
A collection $Ob(Tw_{\Z}(\A))$, $Tw_{\Z}(\alpha,\beta)$ and 
$m_n^{Tw_{\Z}(\A)}$ given above determines a $\Z$-graded $A_\infty$-category 
with $m_n\ne 0$ only if $n=1,2$, i.e.,
a differential graded $($DG$)$ category in the usual sense. 
We denote it by $Tw_{\Z}(\A)$.
\qed
\end{prop}
\begin{defn}\label{defn:zderived}
Let $\A$ be a unital $\Q$-graded $A_\infty$-category of index $h$ 
with $m_n^\A=0$, $n\ge 3$.
We call the cohomology category of $Tw_{\Z}(\A)$ 
the $\Z$-equivariant bounded derived category of $\A$ and 
denote by $D^b_{\Z}(\A)$.
More precisely, the set of objects is given by 
\begin{equation}
Ob(D^b_{\Z}(\A)):=Ob(Tw_{\Z}(\A)),
\end{equation}
and the set of homomorphisms is given by
\begin{align}
{\rm Hom}_{D^b_{\Z}(\A)}(\alpha,\beta):=&
{\rm Ker}(m_1^{Tw_{\Z}(\A)}:Tw_{\Z}(\A)^0(\alpha,\beta)\to 
Tw_{\Z}(\A)^1(\alpha,\beta))\\
&\left/ {\rm Im}(m_1^{Tw_{\Z}(\A)}:Tw_{\Z}(\A)^{-1}(\alpha,\beta)\to 
Tw_{\Z}(\A)^0(\alpha,\beta))\right. .
\end{align}
\end{defn}
Let $T\in Tw_{\Z}(\A)^0(\alpha,\beta)$ be a $m_1^{Tw_{\Z}(\A)}$-closed 
homomorphism.
As in the case for $D^b(\A)$, we define a mapping cone $C(T)$ as an object
\begin{equation}
C(T):=(a[1]\oplus b; Q_{C(T)}),\quad Q_{C(T)}:=
\begin{pmatrix}
Q_{a[1]} & 0\\
T & Q_b
\end{pmatrix}
.
\end{equation}
We define an exact triangle in the category $D^b_{\Z}(\A)$
as a triangle of the form 
\begin{equation}
\alpha\stackrel{T}{\to}\beta \to C(T)\to \alpha[1],
\end{equation}
for some $T\in {\rm Ker}(m_1^{Tw_{\Z}(\A)}:Tw_{\Z}(\A)^0(\alpha,\beta)\to 
Tw_{\Z}(\A)^1(\alpha,\beta))$.
\begin{thm}
The category $D^b_{\Z}(\A)$ endowed with a shift functor 
$[1]$ and the class of exact triangles defined above
becomes a triangulated category.
\end{thm}
\begin{pf}
As in the case for $D^b(\A)$, the proof is essentiality 
the same as the known results in the usual situation.
\qed\end{pf}
\begin{rem}
The twice of the shift functor $[2]$ is not isomorphic 
to the identity functor in $D^b_{\Z}(\A)$.
\end{rem}
Consider the differential graded (DG) functor ${\rm Tot}$ 
as in \cite{bk:1} 
\begin{equation}
{\rm Tot}:Tw_\Z(Tw_{\Z}(\A))\to Tw_{\Z}(\A),\quad 
(\bigoplus_{i=1}^k(a_i;Q_{a_i});T)\mapsto (\bigoplus_{i=1}^k a_i;
Q+T),
\end{equation}
where 
\begin{equation} 
Q:=
\begin{pmatrix}
Q_{a_1}    & 0       & 0 \\
    0      & \ddots  & 0 \\
    0      & 0       & Q_{a_k} &
\end{pmatrix},
\end{equation}
and $T\in Tw_{\Z}(\A)(\oplus_{i=1}^k(a_i;Q_{a_i}),
\oplus_{i=1}^k(a_i;Q_{a_i}))$ $(\subset\widetilde{\A}(\oplus_{i=1}^k a_i,
\oplus_{i=1}^ka_i)$) satisfies 
$$
m_1^{Tw_{\Z}(\A)}(T)+m_2^{Tw_{\Z}(\A)}(T^{\otimes 2})=0.
$$
It is well-defined since 
\begin{align*}
&m_0^{\widetilde{\A}}(1)+m_1^{\widetilde{\A}}(Q+T)
+m_2^{\widetilde{\A}}((Q+T)^{\otimes 2})\\
=&m_0^{\widetilde{\A}}(1)
+\sum_{i=1}^k m_1^{\widetilde{\A}}(Q_{a_i})+
\sum_{i=1}^k m_2^{\widetilde{\A}}(Q_{a_i}^{\otimes 2})
+m_1^{\widetilde{\A}}(T)+
m_2^{\widetilde{\A}}(T\otimes Q)
+m_2^{\widetilde{\A}}(Q\otimes T)
+m_2^{\widetilde{\A}}(T^2)\\
=&m_0^{\widetilde{\A}}(1)
+\sum_{i=1}^k m_1^{\widetilde{\A}}(Q_{a_i})+
\sum_{i=1}^k m_2^{\widetilde{\A}}(Q_{a_i}^{\otimes 2})
+m_1^{Tw_{\Z}(\A)}(T)+m_2^{Tw_{\Z}(\A)}(T^{\otimes 2}).
\end{align*}
Now the following statement is easily shown as in \cite{bk:1} where
they consider the case when $\A$ is a DG category, i.e., the case when
$m_0^\A=0$ and $h=1$ in our terminology.

\begin{prop}
${\rm Tot}$ is an equivalence of DG categories. 
\qed\end{prop}
Indeed, by definition of the twisted complexes and the differential on them,
we see that $\left(Tw_{\Z}(Tw_{\Z}(\A))(A, B), 
m_1^{Tw_{\Z}(Tw_{\Z}(\A))}\right)$ and 
$\left(Tw_{\Z}(\A)({\rm Tot}(A),{\rm Tot}(B)),m_1^{Tw_{\Z}(\A)}\right)$ 
are the same as complexes. 
\begin{cor}
$D^b_{\Z}(\A)$ is an enhanced triangulated category in the sense of 
Bondal-Kapranov \cite{bk:1}. 
\qed\end{cor}

Let us consider our category $D^b_\Z(\A_f)$ a little bit in detail.

\begin{exmp}
Let $\alpha:=(a=a_+\oplus a_-[-1];Q_a)$ and 
$\beta:=(b=b_+\oplus b_-[-1]:Q_b)$ be objects of $Tw_\Z(\A_f)$.
Then the space of homomorphisms is of the following form:
$$
\Phi\in Tw^q_\Z(\alpha,\beta), \quad q\in 2\Z \Leftrightarrow \Phi=
\begin{pmatrix}
\Phi_{++} & 0\\
0 & \Phi_{--}
\end{pmatrix},
\quad \Phi_{\pm\pm}\in \widetilde{\A_f}^q(a_\pm,b_\pm)_+,
$$
$$
\Phi\in Tw^q_\Z(\alpha,\beta),\quad q-1\in 2\Z\Leftrightarrow \Phi=
\begin{pmatrix}
0 & \Phi_{-+}\\
\Phi_{+-} & 0
\end{pmatrix},
\quad \Phi_{\pm\mp}\in \widetilde{\A_f}^{q\mp 1}(a_\pm,b_\mp)_+,
$$
and the coboundary operator $m_1^{Tw(\A_f)}$ becomes the differential
of the usual form
$$
Q_b\Phi-(-1)^q\Phi Q_a,\quad 
\Phi\in Tw^q_\Z(\alpha,\beta).
$$
Note that if $\Phi\in \widetilde{\A_f}^q(\alpha,\beta)$, then 
\begin{equation}\label{eq:euler}
E\Phi -R_\beta\Phi+\Phi R_\alpha=q\Phi,
\end{equation}
where we put 
$$
R_\alpha:={\rm diag}(\frac{2k_1}{h},\dots,\frac{2k_m}{h},
\frac{2l_1}{h}-1,\dots,\frac{2l_m}{h}-1),\quad 
a=\bigoplus_{i=1}^ma\{\frac{2k_i}{h}\}\oplus 
\bigoplus_{i=1}^ma\{\frac{2l_i}{h}\}[-1],
$$
and 
$$
R_\beta:={\rm diag}(\frac{2k'_1}{h},\dots,\frac{2k'_{m'}}{h},
\frac{2l'_1}{h}-1,\dots,\frac{2l'_{m'}}{h}-1),\quad 
b=\bigoplus_{i=1}^{m'}a\{\frac{2k'_i}{h}\}\oplus 
\bigoplus_{i=1}^{m'}a\{\frac{2l'_i}{h}\}[-1].
$$
By integrating the equation \eqref{eq:euler},
we get for $\lambda\in\C$,
\begin{equation}\label{eq:inteuler}
e^{-\lambda R_\beta}\Phi(e^{\lambda\frac{2a_1}{h}}x_1,\dots,
e^{\lambda\frac{2a_n}{h}}x_n)e^{\lambda R_\alpha}=
e^{q\lambda}\Phi(x_1,\dots,x_n).
\end{equation}
This is the analogue of the homogeneity condition discussed in \cite{hw:1}.
Consider the $\Z$-action defined by 
$x_i\mapsto \exp{(2\pi\sqrt{-1}p\cdot a_i/h)}\cdot x_i$, $p\in\Z$. 
It is clear that $f$ is invariant under this $\Z$-action. 
Note also that $\{\frac{2k}{h}\}$, $k\in\Z$ can be considered as the 
irreducible representations of $\Z$.
For a graded twisted complex $\alpha:=(a;Q_a)$, put 
$$
S_\alpha:={\rm diag}(\frac{2k_1}{h},\dots,\frac{2k_{m}}{h},
\frac{2l_1}{h},\dots,\frac{2l_{m}}{h}),\quad 
a=\bigoplus_{i=1}^{m}a\{\frac{2k_i}{h}\}\oplus 
\bigoplus_{i=1}^{m}a\{\frac{2l_i}{h}\}[-1].
$$
Since $Q_a\in \widetilde{\A}^{1}(a,a)_-$, the similar 
equation as \eqref{eq:inteuler} shows that 
$Q_a$ is equivariant with respect to the $\Z$-action, i.e., we have 
\begin{equation}
e^{-\pi\sqrt{-1} S_\alpha}Q_a(e^{\frac{2\pi\sqrt{-1}a_1}{h}}x_1,\dots,
e^{\frac{2\pi\sqrt{-1} a_n}{h}}x_n)
e^{\pi\sqrt{-1}S_\alpha}=Q_a(x_1,\dots,x_n).
\end{equation}
One can show that there is also the $\Z$-action on 
the space of homomorphisms by \eqref{eq:inteuler}.
For $\Phi_\pm\in \widetilde{\A_f}^q(\alpha,\beta)_\pm$, we have
\begin{equation}
e^{-\pi\sqrt{-1} S_\beta}\Phi_\pm(e^{\frac{2\pi\sqrt{-1}a_1}{h}}x_1,\dots,
e^{\frac{2\pi\sqrt{-1} a_n}{h}}x_n)
e^{\pi\sqrt{-1}S_\alpha}=\pm e^{\phi\sqrt{-1}q}\Phi_\pm(x_1,\dots,x_n).
\end{equation}
Therefore, if $\Phi$ is even $($odd$)$, then $\Phi$ is $\Z$-invariant 
if and only if $q\in 2\Z$ $(q-1\in 2\Z)$.
These facts lead us to our definition of $\Z$-equivariant derived category 
$D^b_{\Z}(\A_f)$ of $\A_f$.
Note that the above $\Z$-action on $\A_f$ factors through $\Z/h\Z$. 
The category whose set of objects is the set of $\Z/h\Z$-equivariant 
matrix factorizations and the space of morphisms is $\Z/h\Z$-invariant 
homomorphisms between matrix factorizations are called in physics 
the category of D-branes in Landau-Ginzburg $(\Z/h\Z$-$)$orbifolds
$($see for example \cite{hw:1}$)$
We can construct it by considering the $\Z/h\Z$-equivariant version of 
$D^b(\A_f)$. 
Indeed, we can show that it is equivalent to $D^b_{\Z}(\A_f)/[2]$.
In order to recover the $\Z$-grading by the shift functor, 
we introduced here the translation $\{2/h\}$ and defined a new 
category $D^b_{\Z}(\A_f)$. 
\end{exmp}
\section{$D^b_{\Z}(\A_f)$ and representations of Dynkin quivers}
The following is our main theorem in this paper.
\begin{thm}\label{thm:quiver}
Let us put $f(x):=x^h\in\C[x]$ for $h\ge 2$ and 
consider the unital $\Q$-graded 
$A_\infty$-category $\A_f$ of index $h$.
Then we have the following equivalence of triangulated categories
\begin{equation}
D^b_{\Z}(\A_f)\simeq D^b({\rm mod-}B),
\end{equation}
where $B$ is the path algebra of the following Dynkin quiver of type $A_{h-1}:$
\begin{equation}\label{eq:anquiver}
\bullet_1\to\bullet_2\to \dots \to\bullet_{h-2} \to\bullet_{h-1},
\end{equation}
$($the algebra of upper triangular matrices over $k)$, and 
$D^b({\rm mod-}B)$ is the bounded derived category of finitely generated 
right $B$-modules.
\end{thm}
\begin{rem}
The above equivalence for $h=2$, 
$D^b_{\Z}(\A_{x^2})\simeq D^b({\rm mod-}\C)$, gives the simplest example 
of Kn\"{o}rrer's periodicity.
\end{rem}
\begin{pf}
It is not difficult to see that our category $D^b_{\Z}(\A_f)$ is 
a Krull-Schmidt category, the spaces of homomorphisms are finite dimensional 
and the endomorphism rings of indecomposable objects are local rings.
See, for example, section 5 of \cite{kr:1} for the proof of 
the general $f$ which defines an isolated singularity.
Therefore, we first 
study the set of isomorphism classes of indecomposable objects.
We use the fact 
that the Auslander-Reiten quiver of the category $D^b(\A_f)$ 
of matrix factorizations for $f$ is given by 
\begin{equation}\label{eq:ar1}
[Q_1]\rightleftharpoons [Q_2]\rightleftharpoons \dots 
\rightleftharpoons [Q_{h-2}] \rightleftharpoons [Q_{h-1}],
\end{equation}
where 
$$
Q_l=
\begin{pmatrix}
0 & x^{h-l}\\
x^l & 0
\end{pmatrix}
,\quad l=1,\dots,h-1,i\in\Z,
$$
and each morphism corresponding to the arrow from left to right is given 
by ${\rm diag}(1,x)$ and the one from right to left is
given by ${\rm diag}(x,1)$.
See \cite{ar:1} and also \cite{o:1}. 
Hence we have the following.
\begin{lem}
The set of isomorphism classes of all indecomposable objects of 
$D^b_{\Z}(\A_f)$ is given by 
\begin{equation}
\left\{[M_{l,i}],\quad l=1,\dots,h-1,i\in\Z \right\},
\end{equation}
where
\begin{equation}
M_{l,i}:=\left(a\{\frac{2i}{h}\}\oplus a\{\frac{2(l+i)}{h}\}[-1];
Q_l \right).
\end{equation}
We also have 
\begin{equation}
\begin{pmatrix}
1 & 0\\
0 & x
\end{pmatrix}
\in {\rm Hom}_{D^b_{\Z}(\A_f)}(M_{l,i},M_{l+1,i}),\quad 
\begin{pmatrix}
x & 0\\
0 & 1
\end{pmatrix}
\in {\rm Hom}_{D^b_{\Z}(\A_f)}(M_{l,i},M_{l-1,i+1}),
\end{equation}
and hence 
\begin{equation}
{\rm Hom}_{D^b_{\Z}(\A_f)}(M_{k,i},M_{l,j})\ne 0,\quad \text{only if }
\quad k+2i\le l+2j.
\end{equation}
In particular, 
\begin{equation}
{\rm Hom}_{D^b_{\Z}(\A_f)}(M_{k,0},M_{l,0})=
\left\{\begin{array}{l}
\C, \quad {\rm if}\ k\le l, \\
0, \quad {\rm if}\ k > l.
\end{array} \right.
\end{equation}
\begin{pf}
One can easily show by direct computations.
\qed
\end{pf}
\begin{rem}
Note that $M_{l,i}[1]\simeq M_{h-l,l+i}$.
\end{rem}
\end{lem}
Serre duality holds in our category $D^b_{\Z}(\A_f)$.
\begin{lem}
There are isomorphisms as $\C$-vector spaces
\begin{equation}
{\rm Hom}_{D^b_{\Z}(\A_f)}(M_{k,i},M_{l,j})\simeq 
{\rm Hom}_{D^b_{\Z}(\A_f)}(M_{l,j},M_{k,i-1}[1])^*,
\quad \text{for all} \quad 
1\le k,l\le h-1,\quad i,j\in\Z.
\end{equation}
\end{lem}
\begin{pf}
There is a trace map \cite{kl:2}
\begin{equation}
Tr_{k}:\widetilde{\A_f}^{1-\frac{2}{h}}(M_{k,i},M_{k,i})_-
\to \C,\quad \Phi\mapsto
\frac{1}{h-1}{\rm Res}\left[\frac{Str(d Q_k\cdot \Phi)}{\frac{\p f}{\p x}}
\right],
\end{equation}
where 
$$
Str(d Q_k\cdot \Phi):=\left[(h-k)x^{h-k-1}\Phi_{+-}-kx^{k-1}\Phi_{-+}\right]dx,
\quad \Phi=
\begin{pmatrix}
0 & \Phi_{-+}\\
\Phi_{+-} & 0
\end{pmatrix}.
$$
$Tr_k(\Phi)=0$ if $\Phi=Q_k\Psi+\Psi Q_k$ for 
some $\Psi$ since $dQ_k\cdot Q_k+Q_k\cdot dQ_k=df\cdot 1_{2\times 2}$.
Note that 
$$
\Phi_k:=\begin{pmatrix}
0 & -x^{h-k-1}\\
x^{k-1} & 0
\end{pmatrix}
\in \widetilde{\A_f}^{1-\frac{2}{h}}(M_{k,i},M_{k,i})_-,\quad
Tr_{k}(\Phi_k)=1,\quad Q_k\Phi_k+\Phi_k Q_k=0.
$$
Therefore, under the isomorphism given by 
$$
\widetilde{\A_f}^{1-\frac{2}{h}}(M_{k,i},M_{k,i})_-\simeq
\widetilde{\A_f}^{0}(M_{k,i},M_{k,i-1}[1])_+,\quad
\begin{pmatrix}
0 & \Phi_{-+}\\
\Phi_{+-} & 0
\end{pmatrix}\mapsto
\begin{pmatrix}
\Phi_{+-}& 0\\
0 & -\Phi_{-+}
\end{pmatrix},
$$
$\Phi_k$ determines an element of 
${\rm Hom}_{D^b_{\Z}(\A_f)}(M_{l,j},M_{k,i-1}[1])$.
Moreover, from the knowledge of the Auslander-Reiten quiver 
\eqref{eq:ar1} of $D^b(\A_f)$, 
we see that ${\rm Hom}_{D^b_{\Z}(\A_f)}(M_{l,j},M_{k,i-1}[1])\simeq \C\Phi_k$.

It is not difficult to see that the following pairings
$$
\widetilde{\A_f}^{\frac{2m}{h}}(M_{k,i},M_{l,j})_+\otimes 
\widetilde{\A_f}^{1-\frac{2}{h}-\frac{2m}{h}}(M_{l,j},M_{k,i})_-\to
\widetilde{\A_f}^{1-\frac{2}{h}-\frac{2m}{h}}(M_{k,i},M_{k,i})_-
\stackrel{Tr_k}{\to} \C,
\quad m\in\Z,
$$
induce the perfect pairings
$$
{\rm Hom}_{D^b_{\Z}(\A_f)}(M_{k,i},M_{l,j})\otimes 
{\rm Hom}_{D^b_{\Z}(\A_f)}(M_{l,j},M_{k,i-1}[1])\to
{\rm Hom}_{D^b_{\Z}(\A_f)}(M_{k,i},M_{k,i-1}[1])
\stackrel{Tr_k}{\simeq}\C.
$$
\qed\end{pf}
\begin{rem}
$S:=\{\frac{-2}{h}\}\circ [1]$ is the Serre functor on $D^b_{\Z}(\A_f)$.
In particular, we have $S^h=[h-2]$.
Therefore $D^b_{\Z}(\A_f)$ is a fractional noncommutative 
Calabi-Yau manifold of dimension $1-2/h$ in the sense of \cite{s:1}.
\end{rem}
Combining the above Serre duality and the data of the 
Auslander-Reiten quiver \eqref{eq:ar1} of $D^b(\A_f)$,
we see that there are no higher extensions among $\{M_{l,0}\}$.
\begin{cor}
For $m\ne 0$, we have 
$$
{\rm Hom}_{D^b_{\Z}(\A)}(M_{k,0},M_{l,0}[m])=0, \quad \text{for all}
\quad i,j=1,\dots, h-1.
$$
\qed\end{cor}
\begin{cor}
$D^b({\rm mod-}B)$ is a full triangulated subcategory of $D^b_{\Z}(\A_f)$.
\end{cor}
\begin{pf}
Use the fact that $(M_{1,0},\dots,M_{h-1,0})$ is a strongly exceptional 
collection and 
$$
B\simeq \bigoplus_{i,j=1}^{h-1}{\rm Hom}_{D^b_{\Z}(\A)}(M_{k,0},M_{l,0}).
$$
Since $D^b_{\Z}(\A_f)$ is an enhanced triangulated category, 
we can apply the theorem by Bondal-Kapranov $($\cite{bk:1} Theorem 1$)$.
\qed\end{pf}
Note that the number of indecomposable objects of $D^b_{\Z}(\A_f)/[2]$ 
is 
$$
\#\left\{[M_{l,i}]~|~l=1,\dots,h-1,\quad i\in\Z/h\Z \right\}=(h-1)\cdot h,
$$
which is the number of roots for the root system $A_{h-1}$.
This number coincides with the number of indecomposable objects 
of $D^b({\rm mod-}B)/[2]$ by Gabriel's theorem \cite{g:1}.
Therefore, $D^b({\rm mod-}B)/[2]\simeq D^b_{\Z}(\A_f)/[2]$.
This proves Theorem \ref{thm:quiver}.
\qed\end{pf}
\begin{rem}
The similar proof can be applied for $D_n$ and $E_6,E_7,E_8$ cases since 
the heart of our proof is to use 
the Auslander-Reiten quivers of $D^b(\A_f)$, 
the fact that any matrix factorization over ADE singularities is gradable, 
the Serre duality and the theorem by Gabriel on the number of indecomposables.
They are well-known or can be shown by direct calculations 
with explicit presentations of matrix factorizations.
We shall discuss this in detail in the next paper \cite{kst:1}.
\end{rem}
\section{Stability condition on $D^b_{\Z}(\A_f)$}
In this section, we will briefly discuss on a stability condition 
on $D^b_{\Z}(\A_f)$.
\begin{defn}\label{defn:phase}
Let $\alpha:=(\oplus_{i=1}^n a\{\frac{2k_i}{h}\}\oplus \oplus_{i=1}^n 
a\{\frac{2l_i}{h}\}[-1] ;Q_a)$ be an object of $D^b_{\Z}(\A_f)$ 
such that $Q_a$ is reduced, i.e., each matrix element of $Q_a$ is 
in the maximal ideal generated by $(x_1,\dots,x_n)$.
Then we call the real number
\begin{equation}
\phi_{\alpha}:=\frac{1}{2n}{\rm Tr} S_a -\frac{1}{2}, \quad 
S_a:={\rm diag}(\frac{2k_1}{h},\dots,\frac{2k_n}{h},
\frac{2l_1}{h},\dots,\frac{2l_n}{h})
\end{equation}
phase of the object $\alpha$.
\end{defn}
\begin{exmp}
Let $f:=x^{n+1}$ and consider the objects 
\begin{equation}
M_{l,i}:=\left(a\{\frac{2i}{h}\}\oplus a\{\frac{2(l+i)}{h}\}[-1];
\begin{pmatrix}
0 & x^{h-l}\\
x^l & 0
\end{pmatrix}
 \right).
\end{equation}
Then 
$$
\phi_{M_{l,i}}=\frac{l+2i}{h}-\frac{1}{2}.
$$
\end{exmp}
\begin{defn}\label{defn:centralcharge}
Let $\omega:=\exp{2\pi\sqrt{-1}/h}$.
For $\alpha=(\oplus_{i=1}^n a\{\frac{2k_i}{h}\}\oplus
\oplus_{i=1}^n a\{\frac{2l_i}{h}\}[-1] ;Q_a)$, 
we define a $\C$-linear map $Z_\omega:K_0(D^b_{\Z}(\A_f))\to\C$ as
follows$:$
\begin{equation}
Z_\omega(\alpha):=\sum_{i=1}^n (\omega^{k_i}-\omega^{l_i}).
\end{equation}
\end{defn}
By the above example, we see that $\phi_{M_{l,i}}$ is the 
phase of $Z_\omega([M_{l,i}])$:
\begin{prop}
\begin{equation}
Z_\omega([M_{l,i}])=2\sin(\frac{l}{h}\pi)\cdot e^{\pi\sqrt{-1}\phi_{M_{l,i}}}.
\end{equation}
\qed\end{prop}
Since we know that all indecomposable objects in $D^b_{\Z}(\A_f)$ 
for $f=x^{n+1}$ have definite phases, we can define a stability 
condition on $D^b_{\Z}(\A_f)$.
\begin{thm}\label{thm:stab}
Let $f:=x^{n+1}$ and $P(\phi)$ be the full additive subcategory of 
$D^b_{\Z}(\A_f)$ whose objects have phase $\phi\in\R$.
Then $P(\phi)$ and $Z_\omega$ define a stability condition on 
$D^b_{\Z}(\A_f)$ in the sense of Bridgeland \cite{b:1}.
More precisely, $P(\phi)$ and $Z_\omega$ satisfy 
the following properties$:$
\begin{enumerate}
\item if $M\in P(\phi)$, then $Z_\omega(M)=m(M)\exp(\sqrt{-1}\pi \phi)$ 
for some $m(M)\in\R_{\ge 0}$,
\item for all $\phi\in\R,$ $P(\phi+1)=P(\phi)[1]$,
\item if $\phi_1>\phi_2$ and $M_i\in P(\phi_i)$, then 
${\rm Hom}_{D^b_{\Z}(\A_f)}(M_1,M_2)=0$,
\item for each nonzero object $M\in D^b_{\Z}(\A_f)$, there is a finite
sequence of real numbers 
$$
\phi_1 >\phi_2>\dots >\phi_n
$$
and a collection of exact triangles
$$
M_{i-1}\to M_i\to N_i\to M_{i-1}[1],\quad M_n:=M,\quad M_0:=0
$$
with $N_j\in P(\phi_j)$ for all $j$.
\end{enumerate}
\qed\end{thm}
The space of stability conditions for $D^b_{\Z}(\A_f)$ should be 
isomorphic to the base space of the universal unfolding of $f$ by
the mirror symmetry.
Therefore we expect that there exists a natural Frobenius 
(K.~Saito's flat) structure on the space of stability conditions and 
the stability condition constructed above should correspond to
the origin of the base space of the universal unfolding.
We shall study this in detail elsewhere.

%
%
\end{document}